\documentclass[12pt,bezier]{article}
\usepackage{times}
\usepackage{booktabs}
\usepackage{pifont}
\usepackage{floatrow}
\usepackage{caption}
\usepackage{mathrsfs}
\usepackage{amsmath}
\usepackage{amsfonts,amsthm,amssymb,mathrsfs,bbding}
\usepackage{txfonts}
\usepackage{graphics,multicol}
\usepackage{graphicx}
\usepackage{color}
\usepackage{caption}
\usepackage{cite}
\usepackage{latexsym,bm}
\usepackage{indentfirst}
\usepackage{xcolor}
\usepackage{enumitem}
\usepackage[colorlinks=true,
linkcolor=blue,citecolor=blue,
urlcolor=blue]{hyperref}

\evensidemargin=\oddsidemargin\topmargin=-1.5cm

\newtheorem{lem}{Lemma}[section]
\newtheorem{theorem}{Theorem}[section]

\newtheorem{clm}{Claim}[section]

\theoremstyle{definition}

\addtocounter{section}{0}

\begin{document}
 \title{Asymptotic Order of the Minimum Number of Graph-Saturating Edges I }
\author{Xianglei Chen$^{a}$, Changxiang He$^{a}$\footnote{Corresponding author:cxhe@usst.edu.cn}\; 
\let\thefootnote\relax\footnote{This research is supported by the National Natural Science Foundation of China (Nos. 12271362, 12671416).}\\
{\small $^{a}\,$School of Mathematics, Faculty of Basic Sciences, University of Shanghai for Science and Technology,} \\ 
{\small Shanghai 200093, China} }

\date{}
\maketitle {\flushleft\large\bf Abstract}
Given an $H$-free graph $G$, a non-edge $e$ of $G$ is called  $H$-saturating if $G+e$ contains a copy of $H$. The function $f_H(n,m)$ denotes   the minimum number of $H$-saturating non-edges over  all $n$-vertex $H$-free graphs with exactly $m$ edges.      In this paper, we fully  characterize the asymptotic order  of $f_H(n,\operatorname{ex}(n,R)+1)$  for any  fixed graph $H$ with $\chi(H)\ge \chi(R)$, where $R$ is any fixed edge-critical graph, and resolve the conjecture proposed by  He et al. (  {\it J. Combin. Theory Ser. B} 159 (2023), 42-60 ). 
\begin{flushleft}
\textbf{Keywords:} Saturating edge; Tur\'an number; doubly-edge-$r$-critical graph 
\end{flushleft}
\textbf{AMS Classification:} 05C35\ \ 05C15

\section{Introduction}

Extremal graph theory has its origins in the classical work of Mantel \cite{Mantel},
  and Turán-type problem have since become a classic problem in extremal graph theory. Given a graph $H$,  a graph is said to be \emph{$H$-free } if it contains no copy of $H$ as a subgraph.  The \emph{ Turán number} of $H$, denoted by $\operatorname{ex}(n,H)$, is defined as the maximum number of edges in an $n$-vertex $H$-free graph.   Let $G$ be an $H$-free graph and $e$ be a \emph{ non-edge} of $G$ (i.e., an edge in the complement of $G$). We say that $e$ is an \emph{$H$-saturating edge} of $G$ if the graph $G + e$ contains a copy of $H$. Indeed, an $H$-free graph $G$ is maximal if and only if every non-edge of $G$ is an $H$-saturating edge. These two problems are closely interconnected.

For an $H$-free graph $G$, we define $f_H(n, m)$ as the   minimum number of  \emph{$H$-saturating edges} among all  $n$-vertex $H$-free graphs with $m$ edges.  
Let $\chi(H)$ denote   the chromatic number of $H$. A graph $H$ is  \emph{ edge-critical} if it has an edge $e$ such that $\chi(H-e)= \chi(H)-1$, and such an edge is a \emph{ critical edge} of $H$.
A pair of distinct edges $\{e_1,e_2\}$ of a graph $H$ is called a \emph{$\chi$-critical two-edge set} if $\chi(H-e_1-e_2)=\chi$ and $\chi(H-e_i)>\chi$ for $i=1,2$. Such a set is called a $\chi$ critical \emph{matching} if $e_1$ and $e_2$ are vertex-disjoint, and a $\chi$-critical \emph{cherry} if they share a common endpoint.
Following Jiang and Pikhurko \cite{Jiang}, we call a graph $H$   \emph{ doubly edge-$r$-critical} if $\chi(H-e)\geq r+1$ for every $e\in E(H)$ and there exist distinct edges $e_1,\ e_2$ such that $\chi(H-e_1-e_2)= r$. When  $\chi(H)=r+2$, the first condition is automatic, and the two witnessing edges are necessarily vertex-disjoint.  

Let $\chi(H)=\chi$. Then it follows that $f_H(n, m) = 0$ for all integers $m \leq  ex(n, K_{\chi-1})$. In \cite{2023}, the authors raise the following conjecture: 
$$f_H \left(n,  \operatorname{ex}(n,K_{\chi-1})+1\right) =\Theta (n^2)$$  for every doubly edge $(\chi-2)$-critical graph $H$.  And  they also prose the following problem: whether  $f_H \left(n,  \operatorname{ex}(n,K_{\chi-1})+1\right)$ is determinable for every doubly edge $(\chi-2)$-critical graph $H$ satisfying $\chi(H)\geq 4$.

In this paper, we fully resolve the conjecture. In the process of proving this conjecture, we establish bounds for the coefficient of the leading term of $f_H \left(n, \operatorname{ex}(n,K_{\chi-1})+1\right)$, and  partially answer the second question.

\begin{theorem} \label{high}
Let $\chi\geq 2$ be an integer, let  $H$ be a fixed graph with $\chi(H)>\chi$ and $\chi(H-e)>\chi$, for every $e\in E(H)$. Let $R$ be an arbitrary fixed edge-critical graph satisfying $\chi(R)=\chi+1$. Then, for all sufficiently large integer $n$, 
exactly one of the following holds.
\begin{enumerate}[label=(\arabic*) ] 
\item  If $H$ has a $\chi$-critical two-edge set, then exactly one of the following holds:
\begin{enumerate}[label=(\roman*)]
    \item If $H$ has a $\chi$-critical matching, then
    $$f_{H}\left(n,\operatorname{ex}(n,R)+1\right)=\Theta (n^2).$$
    \item If $H$ has no $\chi$-critical matching, then
$$f_{H}\left(n,\operatorname{ex}(n,R)+1\right)=\left(\frac{2}{\chi}+o(1)\right)n.$$
\end{enumerate}
\item If $H$ has no $\chi$-critical two-edge set, then 
$$f_{H}\left(n,\operatorname{ex}(n,R)+1\right)=0.$$
\end{enumerate}
\end{theorem}



\section{Some useful lemmas}

In this section, we present several useful lemmas that will be used in the subsequent proofs.
The following classical result of Simonovits determines the exact Tur\'an number of edge-critical graphs. 

\begin{lem}[\cite{Simo}]\label{criex}
    If a graph $F$ is edge-critical and $\chi=\chi(F)$, then $\operatorname{ex}(n,F)=t_{\chi-1}(n)$ for every sufficiently large $n$.
\end{lem}

We also recall the following Graph Removal Lemma, which is a fundamental result in extremal graph theory. This lemma provides a connection between the number of copies of a fixed graph and the number of edge deletions required to eliminate all such copies.  

\begin{lem}[\cite{alon,furedi}, Graph Removal Lemma]\label{GRL}
    Let $F$ be a fixed graph. If an $n$-vertex graph contains $o(n^{v(F)})$ copies of $F$, then $G$ can be made $F$-free by deleting $o(n^2)$ edges.
\end{lem}

The following clique-reduction result of Erd\H{o}s, Frankl, and R{\"o}dl establishes a connection between the chromatic number of a fixed graph $F$ and the structure of $F$-free graphs.

\begin{lem}[\cite{EFV}]\label{reduction}
    Let $F$ be a fixed graph with $\chi(F)\geq 3$. Then every $n$-vertex $F$-free graph can be made by $K_{\chi(F)}$-free by deleting $o(n^2)$ edges.
\end{lem}

\section{Proof of Theorem \ref{high}}

This section is devoted to the proof of Theorem \ref{high}. 
Let $\chi\geq 2$ be an integer, and let  $H$ be a fixed graph with $\chi(H)>\chi$ and $\chi(H-e)>\chi$, for every $e\in E(H)$. Let $R$ be an arbitrary fixed edge-critical graph satisfying $\chi(R)=\chi+1$. By Lemma \ref{criex}, $$\operatorname{ex}(n,R)=t_{\chi}(n),$$
for all sufficiently large integer $n$.

We now distinguish two cases according to whether $H$ has a $\chi$-critical two-edge set.



{\bf Case 1} $H$ has a $\chi$-critical two-edge set. 

Fix a $\chi$-critical two-edge set $\{e_1,e_2\}$  of $H$, choosing it to be a matching whenever such a matching exists. For $i\in \{1,2\}$,
$$\chi< \chi(H-e_i)\leq \chi(H-e_1-e_2)+1=\chi+1.$$
Thus
\begin{equation}\label{chi1}
    \chi(H-e_1)=\chi(H-e_2)=\chi+1.
\end{equation}
Fix a proper $\chi$-coloring of $H-e_1-e_2$ with color classes
$$V(H)=P_1\cup P_2\cup\cdots \cup P_{\chi}.$$

In this coloring the endpoints of each $e_i$ lie in the same color class. Indeed, if the endpoints of $e_i$ had different colors, adding $e_i$ to $H-e_1-e_2$ would give a $\chi$-coloring of $H-e_{3-i}$, contradicting (\ref{chi1}). Relabel the classes so that the endpoints of $e_1$ belong to $P_1$, and the endpoints of $e_2$ belong to $P_\ell$, where $\ell\in \{1,2\}$.

Let $G$ be an $n$-vertex $H$-free graph with $e(G)=t_{\chi}(n)+1$, and suppose that $\ f_H(G)=o(n^2).$
We first establish several structural properties that will be used in both subcases below.

Write $$F_H(G)=\{xy\notin E(G): H\subseteq G+xy\},\ f_H(G)=|F_H(G)|.$$
We now construct an auxiliary edge-colored graph. Let $$\Gamma=K_n-F_H(G).$$
Assign a distinct color to each edge of $G$, and color every edge in $E(\Gamma)\backslash E(G)$ with a single additional color, denoted by 0. Then $\Gamma$ contains no rainbow copy of $H$. Indeed, such a copy would contain at most one edge of color 0, and hence would either lie in $G$ or arise by adding a single non-$H$-saturating edge to $G$, both of which are impossible.
Since $G$ has $\left(\frac{1}{2\chi}+o(1)\right)n^2$ non-edges and $f_H(G)=o(n^2)$, there exists an edge of color 0.
Choose one such edge $e_0$, and let $$L=G+e_0.$$
Then $$e(L)=t_{\chi}(n)+2,$$
and all edges of $L$ have pairwise distinct colors.

We next introduce a rooted auxiliary graph $Q$ whose copies in $L$ can be controlled by the $H$-saturating edges of $G$. Let $e_1=ab$. Take two vertex-disjoint copies of $H-e_1$, and construct $Q$  by identifying the two vertices corresponding to $a$ into a single vertex, and similarly identifying the two vertices corresponding to $b$. All other vertices remain distinct, and no additional edges are added. The two identified vertices are called the roots of $Q$. The two copies of $H-e_1$ are edge-disjoint, and 
$$v(Q)=2v(H)-2,\ \chi(Q)=\chi+1.$$
The latter equality follows from $\chi(H-e_1)=\chi+1$, since $Q$ contains $H-e_1$ and admits a proper $(\chi+1)$-coloring obtained by applying the same coloring to both copies of $H-e_1$.

\begin{clm}\label{clm:removal}
    The graph $L$ can be made $Q$-free by deleting $o(n^2)$-edges.
\end{clm}

\begin{proof}
    We first show that the root pair of $Q$ in $L$ belongs to $F_H(G)$.
    Let $x$, $y$ be the roots of a copy of $Q$ in $L$, and let $H_1$, $H_2$ denote its two edge-disjoint copies of $H-e_1$. Suppose, for a contradiction, that $xy\notin F_H(G)$. Then $xy\in E(\Gamma)$. 
    Note that $xy\notin E(L)$. Indeed, if $xy\in E(L)$, $H_1+xy$ is a copy of $H$ in $L$. As all edges of $L$ have distinct colors, this gives a rainbow copy of $H$ in $\Gamma$, contradicting that  $\Gamma$ contains no rainbow copy of $H$.
    Since $G\subseteq L$, the edge $xy$ has color 0.
    Recall that $e_0$ is the unique edge of color 0 in $L$. As $H_1$ and $H_2$ are edge-disjoint, $e_0$ belongs to at most one of them, and it may belong to neither. In either case, at least one of these two copies, say $H_1$, dose not contain $e_0$. Thus all edges of $H_1$ have distinct nonzero colors. Adding $xy$ produces a rainbow copy of $H$ in $\Gamma$, again a contradiction.
    Hence, $xy\in F_H(G)$.

For each unordered root pair in $F_H(G)$, there are at most two assignments of the ordered roots, and the other $v(Q)-2$ vertices at most $n^{v(Q)-2}$ choices.
    It follows that the number $N_Q(L)$ of labeled copies of $Q$ in $L$ satisfies
    $$N_Q(L)\leq 2f_H(G)n^{v(Q)-2}=o(n^{v(Q)}).$$
    Since $Q$ is fixed, Lemma \ref{GRL} directly implies that $L$ can be made $Q$-free by deleting $o(n^2)$ edges.
\end{proof}

Since $\chi(Q)=\chi+1$, we can further apply a clique-reduction result and the stability theorem to obtain the following near-Tur{\'a}n structure of $L$.

\begin{clm}\label{clm:partition}
    There exists a partition $V(L)=A_1\cup A_2\cup\cdots \cup A_{\chi}$
    maximizing the number of edges between distinct parts such that 
    $$\sum_{i=1}^{\chi}e\left(L[A_i]\right)=o(n^2),\ \sum_{1\leq i<j\leq\chi}\left(|A_i||A_j|-e\left(A_i,A_j\right)\right)=o(n^2).$$
    Moreover, $$\sum_{1\leq i<j\leq\chi}\left(|A_i||A_j|-e\left(A_i,A_j\right)\right)\leq \sum_{i=1}^{\chi}e\left(L[A_i]\right)-2,\ |A_i|\geq \frac{n}{2\chi}, \ \text{for all} \ 1\leq i\leq \chi.$$
\end{clm}

\begin{proof}
    By Claim \ref{clm:removal}, $L$ can be made $Q$-free by removing $o(n^2)$ edges. Since $\chi(Q)=\chi+1$, Lemma \ref{reduction} allows us to delete a further $o(n^2)$ edges and obtain a $K_{\chi+1}$-free graph $L_2$. Consequently, $$e(L_2)=t_\chi(n)-o(n^2).$$  
    The standard clique stability theorem (see \cite{Simo}) gives a $\chi$-partition of $V(L_2)$ with $o(n^2)$ internal edges. The same partition has $o(n^2)$ internal edges of $L$, since only $o(n^2)$ edges were deleted. Choose among all $\chi$-partitions of $V(L)$ one maximizing the number of crossing edges, and call its parts $A_1,\ldots,A_\chi$. By maximality,  $$\sum_{i=1}^{\chi}e(L[A_i])=o(n^2).$$ 
    Let $T=K_{A_1,\ldots,A_\chi}$. Using $e(L)=t_\chi(n)+2$ and $e(T)\le t_\chi(n)$, we obtain
     \begin{align*}
     \sum_{1\leq i<j\leq \chi}\bigl(|A_i||A_j|-e_L(A_i,A_j)\bigr)=\sum_{i=1}^{\chi} e(L[A_i])+e(T)-e(L)\le\sum_{i=1}^{\chi} e(L[A_i])-2=o(n^2).
     \end{align*}  
    Also $e(T)\ge e(L)-\sum_i e(L[A_i])=t_\chi(n)-o(n^2)$. 
    Since
    $$e(T)=\frac12\left(n^2-\sum_{i=1}^{\chi}|A_i|^2\right),\ t_\chi(n)=\left(\frac{\chi-1}{2\chi}+o(1)\right)n^2,$$  
     we have $$\sum_{i=1}^{\chi}\left(|A_i|-\frac{n}{\chi}\right)^2=o(n^2).$$
     Thus $$|A_i|=\left(\frac{1}{\chi}+o(1)\right)n$$ for every $i$, giving the asserted lower bound.
\end{proof}

For the remaining counting arguments, let $M$ be the set of missing
crossing pairs relative to $A_1,\ldots,A_\chi$. By Claim
\ref{clm:partition},
\[
|M|\le \sum_{i=1}^{\chi}e(L[A_i])-2,\
|M|=o(n^2),\ |A_i|=\left(\frac{1}{\chi}+o(1)\right)n.
\]
For a vertex $v$, write $\mu(v)$ for the number of pairs in $M$
incident with $v$, and $N_M(Z)$ for all vertices incident with an
$M$-pair whose other endpoint belongs to $Z$. Fix a small constant
\[
0<\eta<\frac{1}{100\chi v(H)},\qquad
W=\{v:\mu(v)\ge\eta n\},\qquad W_i=W\cap A_i.
\]
Since $\sum_v \mu(v)=2|M|$, we have $|W|=o(n)$.

We next establish a rainbow embedding principle. Suppose that the endpoints of $e_1,e_2$ are mapped injectively to two differently colored internal edges of $\Gamma$, respecting the fixed $\chi$-coloring of $H-e_1-e_2$. Assume that all prescribed vertices lie outside $W$ and that every edge of $H-e_1-e_2$ between them is mapped to an edge of $L$. If one of the two internal edges has color $0$, assume further that $e_0$ is not a crossing pair between prescribed vertices.

\begin{clm}\label{clm:rainbow}
Every mapping satisfying the above conditions extends to a rainbow copy of $H$ in $\Gamma$.
\end{clm}
\begin{proof}
Embed the unprescribed vertices of $H-e_1-e_2$ one at a time into the
corresponding parts $A_j\setminus W$.
Every edge of $H-e_1-e_2$ joins distinct color classes, so only missing crossing pairs can obstruct such an embedding.
We may also exclude any
unprescribed endpoint of $e_0$. At each step, every previously embedded
neighbor forbids fewer than $\eta n$ candidates through missing
crossing pairs. Thus the number of available candidates is at least
\[
\left(\frac{1}{\chi}-o(1)\right)n-v(H)\eta n-O_H(1)>0.
\]
This completes an embedding of $H-e_1-e_2$ in $L$. No color appears
twice among its edges, since all edges of $L$ have distinct colors.
If a prescribed internal edge has color $0$, the embedded copy of $H-e_1-e_2$ avoids $e_0$: its unprescribed endpoints were excluded, no prescribed crossing pair equals $e_0$, and edges within a color class do not belong to $H-e_1-e_2$. Hence color $0$ is not repeated when the two prescribed internal edges are added. Hence the resulting copy of $H$ is rainbow, a contradiction to the choice of $\Gamma$.
\end{proof}

\begin{clm}\label{clm:star}
Suppose $W\ne\varnothing$ and, for every $1\le i\le\chi$, at most one
edge of $L[A_i]$ has both endpoints outside $W_i$. Then there exist
$i$ and $x\in A_i$ such that
\begin{equation}\label{eq:star}
d_L(x,A_j)\ge d_L(x,A_i)\ge\frac{\eta}{\chi}n-1
\qquad(1\le j\le\chi).
\end{equation}
\end{clm}
\begin{proof}
Choose $i$ maximizing $e(L[A_i])$ and put
$\Delta=\max_{x\in A_i}d_{L[A_i]}(x)$. Since $W\ne\varnothing$,
$|M|\ge\eta n$. By Claim \ref{clm:partition},
\[
e(L[A_i])\ge \frac1\chi\sum_{j=1}^{\chi}e(L[A_j])
\ge\frac{|M|+2}{\chi}\ge\frac{\eta n+2}{\chi}.
\]
For sufficiently large $n$ this exceeds $1$, so the hypothesis that
at most one internal edge avoids $W_i$ implies $W_i\ne\varnothing$.
By the hypothesis, all but at most one of the edges of $L[A_i]$
meet $W_i$, so
\[
e(L[A_i])\le |W_i|\Delta+1.
\]
All pairs in $M$ are crossing, and therefore no pair of $M$ has
both endpoints in $W_i$. Consequently,
\[
|W_i|\eta n\le\sum_{w\in W_i}\mu(w)
\le |M|<\sum_{j=1}^{\chi}e(L[A_j])
\le\chi e(L[A_i])\le\chi(|W_i|\Delta+1).
\]
It follows that $\Delta\ge\eta n/\chi-1$. Choose $x\in A_i$
with $d_L(x,A_i)=\Delta$. Moving $x$ to any other part cannot
increase the number of crossing edges in a maximum $\chi$-cut;
thus $d_L(x,A_j)\ge d_L(x,A_i)$ for $j\ne i$.
\end{proof}

Whenever a vertex $x$ satisfying (\ref{eq:star}) is available, we shall
use the following completion procedure. Put
$\alpha=\eta/(4\chi)$ and $\tau=\alpha/(4v(H))$. For each $j$,
choose
\[
C_j\subseteq N_L(x)\cap A_j,\qquad |C_j|=\lfloor\alpha n\rfloor,
\]
avoiding $V(e_0)\setminus\{x\}$, and let
\[
U_j=\{v\in C_j:\mu(v)<\tau n\}.
\]
Since $|M|=o(n^2)$, at most $o(n)$ vertices have $M$-degree
at least $\tau n$, whence
\begin{equation}\label{eq:Uj}
|U_j|\ge (\alpha-o(1))n\ge\frac{3\alpha n}{4}
\qquad (1\le j\le\chi).
\end{equation}
A partial embedding of $H-e_1-e_2$ with at most four prescribed
vertices can then be completed greedily, provided the prescribed
vertices respect the color classes, all required edges between them
are in $L$, one prescribed vertex is $x$, and the other prescribed
vertices belong to the appropriate $U_j$. Indeed, $x$ is adjacent
to every $C_j$, each other prescribed vertex has $M$-degree less
than $\tau n$, and at each step the previously embedded neighbors
exclude at most
\[
v(H)\tau n+O_H(1)=\frac{\alpha n}{4}+O_H(1)
\]
candidates. By (\ref{eq:Uj}), at least $3\alpha n/4-\alpha n/4-O_H(1) =\alpha n/2-O_H(1)>0$ candidates remain, so the embedding can be completed. Since the
sets $C_j$ avoid $V(e_0)\setminus\{x\}$, the embedded copy of
$H-e_1-e_2$ uses no edge $e_0$.

We now distinguish two subcases according to whether $H$ has a $\chi$-critical matching.

{\bf Subcase 1.1} $H$ has a $\chi$-critical matching.

Suppose the asserted quadratic lower bound fails. Then the graphs
$G$ above may be chosen with $f_H(G)=o(n^2)$. By our choice of
$e_1,e_2$, they form a $\chi$-critical matching. The next
claim verifies the hypothesis of Claim \ref{clm:star}.

\begin{clm}\label{clm:matching-concentration}
The set $W$ is not empty. Moreover, 
$e(L[A_i\backslash W_i])\le 1$
for every $1\le i\le\chi$.
\end{clm}
\begin{proof}
First suppose the endpoints of $e_1$ and $e_2$ have different
colors in the fixed proper coloring of $H-e_1-e_2$. We claim that
some $\Gamma[A_i]$ contains two differently colored edges. Otherwise,
each $L[A_i]$ has at most one edge, as $L$ contains only one
representative of each color. By Claim \ref{clm:partition},
\[
2\le\sum_i e(L[A_i])\le\chi,
\qquad |M|\le\sum_i e(L[A_i])-2.
\]
If every pair of internal edges lying in different parts had a
missing crossing pair between their endpoints, then, since each
part contains at most one internal edge and each missing crossing
pair can obstruct at most one pair of these edges,
\[
|M|\ge\binom{\sum_i e(L[A_i])}{2}
\ge\sum_i e(L[A_i])-1,
\]
a contradiction.
Therefore two differently colored internal edges in distinct parts
have all crossing pairs between their endpoints in $L$. Since
$|M|=O_H(1)$ in this situation, their endpoints lie outside $W$
for all sufficiently large $n$. Claim \ref{clm:rainbow} gives a
contradiction.

Next, suppose two differently colored edges of $\Gamma[A_i]$ have
all their endpoints outside $W$. Let $Z$ be their set of at most
four endpoints and choose $j\ne i$. The set
\[
U=A_j\setminus\bigl(W\cup N_M(Z)\cup V(e_0)\bigr)
\]
has $|U|\ge(1/\chi-4\eta-o(1))n>n/(2\chi)$ for all large $n$.
Since $f_H(G)=o(n^2)$, there is a pair $uv\in\binom U2$ belonging
to $E(\Gamma)$. At least one of the original two edges has color
different from $uv$. That edge and $uv$ lie in distinct parts and
all crossing pairs between their endpoints belong to $L$.
All prescribed vertices are outside $W$, and $uv$ avoids $V(e_0)$.
Claim \ref{clm:rainbow} again gives a contradiction. The first
paragraph now shows that $W\ne\varnothing$. Moreover, since the
internal edges of $L$ have different colors, at most one edge of
$L[A_i]$ avoids $W_i$, for every $i$.

Now suppose both critical edges have the same color in the fixed
coloring. By Claim \ref{clm:rainbow}, no two vertex-disjoint,
differently colored edges of $\Gamma[A_i]$ have all their endpoints
outside $W$. Since $\sum_i e(L[A_i])\ge2$, there is an internal edge $uv\in L[A_i]$
with a private nonzero color. If $W=\varnothing$, the set
$A_i\setminus(V(uv)\cup V(e_0))$ has linearly many vertices and
contains a pair from $E(\Gamma)$, since $f_H(G)=o(n^2)$.
This pair is disjoint from $uv$ and has a different color,
a contradiction. Hence $W\ne\varnothing$.

More generally, if a nonzero-colored edge $uv\in L[A_i]$ avoids
$W_i$, the set $A_i\setminus(W\cup V(uv)\cup V(e_0))$ still has
linearly many vertices, because $|W|=o(n)$. It contains a pair
from $E(\Gamma)$, giving the same contradiction. Thus every
edge of $L[A_i]$ outside $W_i$ must be $e_0$, and there can be
at most one such edge. This proves the claim in both color
configurations.
\end{proof}

By Claims \ref{clm:matching-concentration} and \ref{clm:star},
there exist $i$ and $x\in A_i$ satisfying (\ref{eq:star}).
Choose $C_j,U_j$ as above. If $e_1$ and $e_2$ have different
forced colors, choose $j\ne i$. As $|U_j|=\Omega(n)$ and
$f_H(G)=o(n^2)$, there is a pair
$uv\in\binom{U_j}{2}\setminus F_H(G)$.
Since $\mu(u),\mu(v)<\tau n$ and $2\tau<3\alpha/4$,
we may choose $y\in U_i$ such that $yu,yv\in E(L)$.
All four crossing pairs between $\{x,y\}$ and $\{u,v\}$
are edges of $L$. Prescribe the endpoints of $e_1$ as $x,y$ and
those of $e_2$ as $u,v$, identifying their color classes with
$A_i,A_j$ respectively. By the preceding greedy completion,
$H-e_1-e_2$ embeds in $L$ without using $e_0$.
Moreover, $xy\in E(G)$ because $xy\in E(L)$ and $xy\ne e_0$.
Thus $G+uv$ contains $H$, contradicting either the $H$-freeness
of $G$ (if $uv\in E(G)$) or $uv\notin F_H(G)$ (otherwise).

If the two critical edges have the same forced color, take
$uv\in\binom{U_i}{2}\setminus F_H(G)$ and
$y\in U_i\setminus\{u,v\}$. The prescribed endpoints
$x,y,u,v$ belong to the same independent color class of
$H-e_1-e_2$. The greedy completion again embeds $H-e_1-e_2$
in $L$ without $e_0$, while $xy\in E(G)$.
This leads to the same contradiction. Therefore no sequence with
$f_H(G)=o(n^2)$ exists, proving
\[
f_H(n,t_\chi(n)+1)=\Omega(n^2).
\]

{\bf Subcase 1.2} $H$ has no $\chi$-critical matching.

Since $H$ has a $\chi$-critical two-edge set but no $\chi$-critical matching, $e_1$, $e_2$ form a $\chi$-critical cherry.
Write $e_1=ab$ and $e_2=ac$. Fix an arbitrary
$0<\varepsilon<2/\chi$, and suppose for a contradiction that
$f_H(G)<(2/\chi-\varepsilon)n$ for arbitrarily large $n$.
The common reduction applies since $f_H(G)=o(n^2)$.
The three vertices $a,b,c$ lie in
one color class of $H-e_1-e_2$.

We first show that $W\neq \emptyset$ and 
$e(L[A_i\backslash W_i])\le 1$
for every $1\le i\le\chi$.
It suffices to show that every non
zero-colored internal edge of $L$ meets $W$.
Suppose, to the contrary, that a nonzero-colored edge $uv\in E(L[A_i])$ has both endpoints
outside $W$. Then
\[
U=A_i\setminus\bigl(W\cup\{u,v\}\cup V(e_0)\bigr)
\]
has $|U|=(1/\chi-o(1))n$.
For any $w\in U$, consider both pairs $uw$ and $vw$. If either
pair, say $uw$, belongs to $E(\Gamma)$, then $uv$ and $uw$
are differently colored adjacent edges of $\Gamma[A_i]$ with
all three endpoints outside $W$. They represent the cherry
$ab,ac$ in the fixed coloring of $H-e_1-e_2$, and
Claim \ref{clm:rainbow} gives a rainbow copy of $H$, a contradiction.
The same argument applies to $vw$, with $v$ as the common endpoint.
Thus both $uw$ and $vw$ belong to $F_H(G)$ for every $w\in U$.
These $2|U|$ non-edges are pairwise distinct. Consequently,
\[
f_H(G)\ge 2|U|=\left(\frac{2}{\chi}-o(1)\right)n
>\left(\frac{2}{\chi}-\varepsilon\right)n,
\]
contrary to our assumption.
Since $\sum_i e(L[A_i])\ge2$ and $e_0$ is the only color-$0$ edge of $L$,
there is a nonzero-colored internal edge. It meets $W$, so
$W\ne\varnothing$. Every internal edge outside $W_i$ must
be $e_0$.

By Claim \ref{clm:star},
there exist $i$ and $x\in A_i$ satisfying (\ref{eq:star}).
Choose $C_j,U_j$ as above and take distinct $y,z\in U_i$.
Prescribe $a,b,c$ as $x,y,z$, respectively. These three
vertices form an independent set in $H-e_1-e_2$, so the common
greedy completion gives an embedding of $H-e_1-e_2$ in $L$
that avoids $e_0$. Both $xy,xz$ belong to $E(G)$, since the
sets $C_j$ avoid $V(e_0)\setminus\{x\}$. Hence $G$ already
contains a copy of $H$, a contradiction. As $\varepsilon>0$ was arbitrary, this proves
\[
\liminf_{n\to\infty}\frac{f_H(n,t_\chi(n)+1)}{n}
\ge\frac{2}{\chi}.
\]

To complete the proof, choose a smallest part $A$ of $T_\chi(n)$,
take distinct vertices $u,v\in A$, and put $G_0=T_\chi(n)+uv$.
Then $e(G_0)=t_\chi(n)+1$. Moreover $G_0$ is $H$-free: a copy
of $H$ in $G_0$ would become $\chi$-colorable upon deleting
at most one of its edges, contrary to the assumptions of
Theorem \ref{high}.

In Subcase 1.1, the bound $f_H(G_0)\le\binom n2=O(n^2)$ combines
with the quadratic lower bound to give conclusion~(1)(i).
In Subcase 1.2, any copy of $H$ in $G_0+xy$ for a saturating
non-edge $xy$ must use both $xy$ and $uv$; otherwise deleting
at most one edge of $H$ would make it $\chi$-colorable.
The two preimage edges form a $\chi$-critical two-edge set.
Since $H$ has no $\chi$-critical matching, they must be adjacent.
Consequently $xy$ shares an endpoint with $uv$ and lies in $A$.
It follows that
\[
f_H(G_0)\le 2(|A|-2)
=2\left\lfloor\frac n\chi\right\rfloor-4=O(n).
\]
Together with the lower bound in Subcase 1.2, this yields
\[
f_H(n,t_\chi(n)+1)=\left(\frac{2}{\chi}+o(1)\right)n,
\]
proving conclusion~(1)(ii).


{\bf Case 2} $H$ has no $\chi$-critical two-edge set.

Let $G_0$ be obtained from  $T_{\chi}(n)$ by adding an edge $uv$ within one part. Then
$$e(G_0)=t_{\chi}(n)+1=\operatorname{ex}(n,R)+1.$$
Note that $G_0$ is $H$-free. 
Indeed, any copy of $H$ in $G_0$ must contain $uv$, since $\chi(H)>\chi$.
Deleting $uv$ from such a copy would imply that $\chi(H-e')\leq \chi$ for some $e'\in E(H)$, contradicting the assumption that $\chi(H-e)>\chi$ for every $e\in E(H)$.
We next show that $G_0$ has no $H$-saturating edges. Suppose, to the contrary, that $xy$ is an $H$-saturating edge of $G_0$. Then $G_0+xy$ contains a copy of $H$, which must contain $xy$. We claim that such a copy must also contain $uv$. Indeed, suppose otherwise. Let $\varphi$ be an embedding of $H$ into $G_0+xy$ whose image contains $xy$ but not $uv$, and let $e\in E(H)$ be the edge mapped to $xy$. Then $\varphi$ induces an embedding of $H-e$ into $T_{\chi}(n)$. Consequently, $\chi(H-e)\leq \chi$, which is impossible since $\chi(H-e)>\chi$.

Let $e_1,e_2\in E(H)$ be the edges corresponding to $xy$ and $uv$, respectively. Since $H-e_1-e_2$ emdeds into $T_{\chi}(n)$, we have $$\chi(H-e_1-e_2)\leq \chi.$$
On the other hand, deleting one edge decreases the chromatic number by at most one, so $$\chi(H-e_1-e_2)\geq \chi(H-e_1)-1\geq \chi.$$
Consequently,
$$\chi(H-e_1-e_2)= \chi,$$
which contradicts that $H$ has no $\chi$-critical two-edge set. Therefore, $f_{H}\left(n,\operatorname{ex}(n,R)+1\right) =0 $.
This completes the proof of Theorem \ref{high}.

 \section{Concluding remarks}
In the quadratic case of Theorem \ref{high}, our proof shows that there exists a constant $c_{H,\chi}>0$ such that
$$c_{H,\chi}\le \liminf_{n\rightarrow \infty}\frac{f_H(n,\operatorname{ex}(n,R)+1)}{n^2}\le \limsup_{n\rightarrow \infty}\frac{f_H(n,\operatorname{ex}(n,R)+1)}{n^2}\le\frac{1}{2\chi}.$$

In fact, in a separate manuscript entitled \emph{ Asymptotic Order of the Minimum Number of Graph-Saturating Edges II}, Xianglei Chen and Changxiang He \cite{CHII} determine the exact asymptotic leading coefficient in the quadratic case and obtain the exact value in the linear case for all sufficiently large $n$.

\section{Declaration of AI use}
During the preparation of this work, the authors used OpenAI ChatGPT to assist in the development of certain mathematical proofs. All AI-assisted arguments were independently verified by authors, who take full responsibility for the correctness of proofs and the content of manuscript.

\end{document}